\documentclass{gtmon_a}
\pdfoutput=1

\proceedingstitle{The Zieschang Gedenkschrift}
\conferencestart{5 September 2007}
\conferenceend{8 September 2007}
\conferencename{Conference in honour of Heiner Zieschang}
\conferencelocation{Toulouse, France}

\editor{Michel Boileau}
\givenname{Michel}
\surname{Boileau}

\editor{Martin Scharlemann}
\givenname{Martin}
\surname{Scharlemann}

\editor{Richard Weidmann}
\givenname{Richard}
\surname{Weidmann}

\title{A Magnus theorem for some one-relator groups}

\author{Oleg Bogopolski}
\givenname{Oleg} \surname{Bogopolski}
\address{{\rm OB:}\qua
Fachbereich Mathematik\\ Universit{\"a}t Dortmund\\\newline
Vogelpothsweg 87\\ Dortmund, 44227\\Germany\vspace{3pt}\\\newline
{\rm OB, KS:}\qua Institute of Mathematics of SBRAS\\\newline
Koptjuga 4\\ Novosibirsk, 630090\\Russia}
\email{groups@math.nsc.ru}
\urladdr{http://www.mathematik.uni-dortmund.de/lsvi/bogopolski/index.html}

\author{Konstantin Sviridov}
\givenname{Konstantin} \surname{Sviridov}
\email{konstantin\_sviridov@academ.org}
\urladdr{}

\volumenumber{14}
\issuenumber{}
\publicationyear{2008}
\papernumber{04}
\startpage{63}
\endpage{73}

\doi{}
\MR{}
\Zbl{}

\arxivreference{} %%% Not in ArXiV


\received{29 April 2006}
\revised{}
\accepted{9 November 2006}
\published{29 April 2008}
\publishedonline{29 April 2008}
\proposed{}
\seconded{}
\corresponding{}
\version{}

\makeatletter
\def\cnewtheorem#1[#2]#3{\newtheorem{#1}{#3}[section]
\expandafter\let\csname c@#1\endcsname\c@thm}
\makeatother

\AtBeginDocument{\let\bar\wbar\let\tilde\wtilde\let\hat\what}

\newtheorem{thm}{Theorem}[section]
\cnewtheorem{prop}[thm]{Proposition}
\cnewtheorem{lem}[thm]{Lemma}
\cnewtheorem{cor}[thm]{Corollary}
\newtheorem{mainthm}{Main Theorem}

\theoremstyle{definition}
\cnewtheorem{example}[thm]{Example}
\cnewtheorem{problem}[thm]{Problem}
\cnewtheorem{conj}[thm]{Conjecture}
\cnewtheorem{defn}[thm]{Definition}
\cnewtheorem{rem}[thm]{Remark}
\cnewtheorem{qs}[thm]{Question}
\makeautorefname{defn}{Definition}
\makeautorefname{mainthm}{Main Theorem}

\begin{document}

\begin{asciiabstract}
We will say that a group G possesses the Magnus property if for
any two elements u,v in G with the same normal closure, u is
conjugate to v or v^{-1}. We prove that some one-relator
groups, including the fundamental groups of closed nonorientable
surfaces of genus g>3 possess this property. The analogous
result for orientable surfaces of any finite genus was obtained 
by the first author [Geometric methods in group theory, Contemp. Math,
372 (2005) 59-69].
\end{asciiabstract}

\begin{htmlabstract}
We will say that a group G possesses the Magnus property if for
any two elements u,v&isin; G with the same normal closure, u is
conjugate to v or v<sup>-1</sup>. We prove that some one-relator
groups, including the fundamental groups of closed nonorientable
surfaces of genus g>3 possess this property. The analogous
result for orientable surfaces of any finite genus was obtained by
the first author [Geometric methods in group theory, Contemp. Math,
372 (2005) 59-69].
\end{htmlabstract}

\begin{abstract}
We will say that a group $G$ possesses the Magnus property if for
any two elements $u,v\in G$ with the same normal closure, $u$ is
conjugate to $v$ or $v^{-1}$. We prove that some one-relator
groups, including the fundamental groups of closed nonorientable
surfaces of genus $g>3$ possess this property. The analogous
result for orientable surfaces of any finite genus was obtained by
the first author \cite{B}.
\end{abstract}

\begin{webabstract}
We will say that a group $G$ possesses the Magnus property if for
any two elements $u,v\in G$ with the same normal closure, $u$ is
conjugate to $v$ or $v^{-1}$. We prove that some one-relator
groups, including the fundamental groups of closed nonorientable
surfaces of genus $g>3$ possess this property. The analogous
result for orientable surfaces of any finite genus was obtained by
the first author [Geometric methods in group theory, Contemp. Math,
372 (2005) 59-69].
\end{webabstract}

\maketitle

\section{Introduction}\label{sec1}

In 1930 W\,Magnus published a very important (for combinatorial
group theory and logic) article where he proved the so-called 
\textit{Freiheitssatz} and the following theorem.

\begin{thm}{\rm \cite{M}}\qua\label{free}
 Let $F$ be a free group and $r,s\in F$. If the
normal closures of $r$ and $s$ coincide, then $r$ is conjugate to
$s$ or $s^{-1}$.
\end{thm}

In \cite{BKZ}, O\,Bogopolski, E\,Kudryavtseva and H\,Zieschang
proved the analogous result for fundamental groups of closed
orientable surfaces in case where $r$ and $s$ are represented by
simple closed curves. The suggested proof was geometrical and used
coverings, intersection number of curves and Brouwer's fixed-point
theorem. However, they were not able to generalize it for arbitrary
elements $r,s$.

Later O\,Bogopolski, using algebraic methods in the spirit of
Magnus, proved the desired result without restrictions on $r,s$.

\begin{thm}{\rm \cite{B}}\qua\label{or}
 \textit{Let $G$ be the fundamental group of a closed
orientable surface and $r,s\in G$. If the normal closures of $r$
and $s$ coincide, then $r$ is conjugate to $s$ or $s^{-1}$.}
\end{thm}

In \cite{H}, Howie proposed another, topological, proof of this
theorem. Both proofs do not work in the nonorientable case. The
main result of the present article is the following theorem.

\begin{mainthm}\label{mainthm} Let $G=\langle a,b,y_1,\dots ,y_e\,|\,
[a,b]uv\rangle$, where $e\geqslant 2$, $u,v$ are nontrivial
reduced words in letters $y_1,\dots ,y_e$, and $u,v$ have no
common letters. Let $r,s\in G$ be two elements with the same
normal closures. Then $r$ is conjugate to $s$ or $s^{-1}$.
\end{mainthm}

It is known that the fundamental group of a closed nonorientable
surface of genus $k\geqslant 3$ has the presentation $\langle
x_1,x_2,\dots ,x_k\,|\, [x_1,x_2]x_3^2\cdot \ldots \cdot
x_k^2\rangle.$ So, we have the following corollary.

\begin{cor}\label{nonor} Let $G$
%$G=\langle a,b,c,d,\dots ,e,f\,|\, [a,b][c,d]\dots [e,f]=1\rangle$
be the fundamental group of a closed nonorientable surface of
genus at least 4, and $r,s\in G$. If the normal closures of $r$
and $s$ coincide, then $r$ is conjugate to $s$ or $s^{-1}$.
\end{cor}

Note that this corollary trivially holds for genus 1 and 2, but we
do not know, whether it holds for genus 3.

We will say that a group $G$ possesses \textit{the Magnus property},
if for any two elements $r,s$ of $G$ with the same normal closures
we have that $r$ is conjugate to $s$ or $s^{-1}$. So, all the
above theorems imply that the fundamental group of any compact
surface, except of the nonorientable surface of genus 3, possesses
the Magnus property.

It was shown in \cite{B} that the Magnus property does not hold
for many one-relator groups, including generalized
Baumslag--Solitar groups, all noncyclic one-relator groups with
torsion, and infinitely many one-relator torsion-free hyperbolic
groups.

Now we discuss some logical aspects concerning this property. It
was noticed in $\cite{B}$ that if two groups $G_1,G_2$ are
elementary equivalent and $G_1$ possesses the Magnus property,
then $G_2$ possesses this property too. In particular, any group,
which is elementary equivalent to a free group or a free abelian
group possesses the Magnus property. This gives another way of
proving of \fullref{or} and \fullref{nonor}. However,
there are groups, which are not even existentially equivalent to a
free group (hence, they are not limit groups), but possess the
Magnus property. The easiest example is the direct product
$F_n\times F_m$ of nontrivial free groups of ranks $n,m$, where $n+m\geqslant
3$. The other example is the following:
$$G=\langle a,b,x_1,\dots ,x_n,y_1,\dots ,y_m\,|\,
[a,b][X,Y]Z^k\rangle,$$ where $k\geqslant 4$, $X,Y$ are words in
the letters $x_1,\dots, x_n$, and $Z$ is a word in the letters
$y_1,\dots ,y_m$, such that $[X,Y]\neq 1$ and $Z\neq 1$ in the
corresponding free groups. This group possesses the Magnus
property by our \fullref{mainthm}, but is not existentially equivalent
to a free group. Indeed, by \cite{CCE}, for any $l>1$ the $l$--th
power of a nontrivial element of a free group can not be
expressed as a product of less than $\upnfrac{l+1}{2}$ commutators.
Thus, the following formula is valid in $G$, but is not valid in
any free group:
$$\exists\, z_1,z_2,z_3,z_4,z \, ( z\neq 1 \wedge [z_1,z_2][z_3,z_4]z^k=1).$$

\textbf{Problems}\qua (1)\qua Does every amalgamated product
$A\ast_{\mathbb{Z}} B$, where $A,B$ are free groups and
$\mathbb{Z}$ is a maximal cyclic subgroup in both factors,
possesses the Magnus property?

(2)\qua Does every limit group possesses the Magnus property?

(3)\qua Does the group $G=\langle a,b,c\,|\, a^2b^2c^2\rangle$
possesses the Magnus property?

(4)\qua Let $A$ and $B$ be torsion free groups which possess the Magnus
property. Does the group $A\ast B $ possesses the Magnus property?
(A positive answer in a partial case can be found in the paper by Edjvet \cite{E}.
Note also, that the Magnus property is closed under direct
products.)

 Some other problems related to the Magnus property are
collected in \cite{B}.
\medskip

The plan of this paper is the following. In \fullref{section 2} we deduce
the \fullref{mainthm} from \fullref{2.1} and prove auxiliary
\fullref{2.2}. In \fullref{section 3} we introduce some technical notions
like the left and the right bases of a subgroup, the width of an
element, a piece of an element, a special element, and prove
auxiliary \fullref{3.2} and \mbox{\fullref{3.4}}. In \fullref{section 4}
we present some quotients as amalgamated products and prove the
crucial \fullref{4.1}. In \fullref{section 5} we prove \fullref{2.1}.

\subsubsection*{Acknowledgements}
The first author was partially supported by INTAS grant N 03-51-3663
and by the grant Complex integration projects of SB RAS N 1.9.

\section{Some reduction}\label{section 2}

First we introduce notation. Let $A$ be a group, $g,h\in A$. The
normal closure of $g$ in $A$ is denoted by $\langle\!\langle
g\rangle\!\rangle_A$ or simply $\langle\!\langle
g\rangle\!\rangle$ if the group is clear from the context.
%We will write $g\sim_A h$ if $g$ is conjugate to $h$ in $A$.
Denote $[g,h]=g^{-1}h^{-1}gh$. Let $X$ be an alphabet, $x\in X$
and $r$ be a word in the alphabet $X\cup X^{-1}$. By $r_x$ we
denote the exponent sum of $x$ in $r$.

%The fundamental group of the nonorientable surface of genus $g>2$
%has the following presentation
%$$G=\langle a,b,c,d,\dots ,f\,|\,[a,b]c^2d^2\cdot \ldots \cdot
%f^2\rangle,$$ where the number of letters in this presentation is
%equal to $g$. We assume that $g\geqslant 4$. Thus, the letters
%$a,b,c,d$ in this presentation are present. This assumption will
%be used to prove Lemma 2.2. Note that for genus 3 this Lemma is
%not valid.

We will deduce the \fullref{mainthm}  from the following proposition.

\begin{prop}\label{2.1} Let $H=\langle x,b,y_1,\dots
,y_e\,|\, [x^k,b]uv\rangle$, where $e\geqslant 2$, $k\neq 0$,
$u,v$ are nontrivial reduced words in $y_1,\dots ,y_e$, and $u,v$
have no common letters. Let $r,s\in H$ be two elements with the
same normal closures and let $r_x=0$. Then $r$ is conjugate to $s$
or $s^{-1}$.
\end{prop}

\begin{proof}[Proof of the \fullref{mainthm}] Let $r,s\in G$ and the normal
closures of $r$ and $s$ coincide.  Suppose that $r_b=0$. In this
case we will use another presentation of $G$:
$$G=\langle a,b,y_1,\dots ,y_e\,|\, [b,a]v^{-1}u^{-1}\rangle.$$
Then the \fullref{mainthm} follows immediately from \fullref{2.1}.
% Consider the automorphism $\varphi$ of $G$
%given by the rule:

%$\varphi : \left\{
%\begin{array}{l}
%a\mapsto b\\
%b\mapsto a\\
%c\mapsto f^{-1}\\
%d\mapsto e^{-1}\\
%\dots\\
%e\mapsto d^{-1}\\
%f\mapsto c^{-1}.
%\end{array}
%\right.$

%Then $(\varphi(r))_a=0$. By previous case $\varphi(r)$ is
%conjugate to $\varphi(s)^{\pm 1}$, hence $r$ is conjugate to
%$s^{\pm 1}$.

Now suppose that $r_b\neq 0$. In this case we can embed naturally
the group $G$ into the group
$$H=G\underset{a=x^{r_b}}{\ast} \langle x\,|\, \rangle,$$
where $x$ is a new letter. Clearly, the normal closures of $r$ and
$s$ in $H$ coincide.  To finish the proof, we need the following claim.

\medskip
\textbf{Claim}\qua{\sl The elements $r$ and $s$ are conjugate in $H$ if and
only if they are conjugate in~$G$.}

\begin{proof}  Suppose that $r=h^{-1}sh$, where $h\in H$. Write
$h=g_1z_1\dots g_nz_ng_{n+1}$, where $z_i\in \{x,x^2,\dots,
x^{|r_b|-1}\}$, $g_i\in G$ and $g_2,\dots, g_n$ are nontrivial
($g_1$ and $g_{n+1}$ may be trivial). We may assume that $n$ is
minimal possible. Suppose that $n\geqslant 1$. From the normal
form we deduce that $g_1^{-1}sg_1\in \langle a\rangle$. Then $z_1$
centralizes this element, that contradicts to the minimality of
$n$. Hence, $n=0$ and $h\in G$. \end{proof}

So, we will work now with the group $H$. This group has the
following presentation:
$$\langle x,b,y_1,\dots ,y_e\,|\, [x^{r_b},b]uv\rangle.$$ Let $\bar{b}=x^{r_a}b$. Using Tietze
transformation we can rewrite this presentation as
$$\langle x,\bar{b},y_1,\dots ,y_e\,|\,
[x^{r_b},\bar{b}]uv\rangle.$$ Writing $r$ in the generators of
this presentation, we have $r_x=r_ar_b-r_br_a=0$. Again, by
\fullref{2.1},
%we have $r\sim s^{\pm 1}$ in $H$, and hence in $G$.
$r$ is conjugate to $s^{\pm 1}$ in $H$ and, hence in $G$. \end{proof}

Let $c_1,\dots ,c_p$ be the letters of the word $u$, and
$d_1,\dots ,d_q$ be the letters of the word $v$. Consider the
following automorphism of $H$:

$\psi \co \left\{
\begin{array}{l}
x\mapsto x\\
b\mapsto bu\\
c_i\mapsto x^{-k}c_ix^k \hspace*{31mm}(i=1,\dots, p)\\
d_j\mapsto x^{-k}u^{-1}x^{k}d_jx^{-k}ux^k \hspace*{1cm}(j=1,\dots,
q).
\end{array}
\right.$

%$\psi : \left\{
%\begin{array}{l}
%x\mapsto x\\
%b\mapsto bc^2\\
%c\mapsto x^{-k}cx^k\\
%d\mapsto x^{-k}c^{-2}x^{k}dx^{-k}c^2x^k\\
%\dots\\
%f\mapsto x^{-k}c^{-2}x^{k}fx^{-k}c^2x^k.
%\end{array}
%\right.$

\begin{lem}\label{2.2} Let $h$ be a nontrivial element of $H$. Then there exists a
natural $n_0$ such that for all $n>n_0$ the element $\psi^n(h)$ is
not conjugate to a power of $b$.
\end{lem}

\begin{proof} Suppose that there exists an $m$ such that
$\psi^m(h)$ is conjugate to a nonzero power of $b$. Then for any
$l>0$ the element $\psi^{m+l}(h)$ is conjugate to a nonzero power
of $b\prod_{i=0}^{l-1} x^{-ik}ux^{ik}$. But any such power is not
conjugate to a power of $b$ since its image in $H/\langle\!\langle
b,x\rangle\!\rangle$ is nontrivial. \end{proof}

\section{Left and right bases}\label{section 3}

%\textit{Proof of Proposition 1.}
Without loss of generality, we may assume that $k>0$. Consider the
homomorphism $ H\rightarrow \mathbb{Z}$, which sends $x$ to 1 and
any other generator of $H$ to 0. Denote its kernel by~$N$. Denote
$g_i=x^{-i}gx^i$ and $Y_i=\{b_i,(y_1)_i,\dots ,(y_e)_i\}$. Using
the Reidemeister--Schreier method we can find the following
presentation of $N$:
$$N=\Big\langle \bigcup_{i\in \mathbb{Z}} Y_i\,|\, b_iu_iv_i=b_{i+k}\, (i\in \mathbb{Z})\Big\rangle.$$
Denote by $w_i$ the word $b_iu_iv_i$. Thus, $w_i=b_{i+k}$ in $N$.
We will use the following three presentations of $N$ depending on
a situation:

(a)\qua $N$ is the free product of the free groups $G_i=\langle
Y_i\,|\, \rangle$\, $(i\in \mathbb{Z})$ with amalgamation, where
$G_i$ and $G_{i+k}$ are amalgamated over the cyclic subgroup
generated by $w_i$ in $G_i$ and by $b_{i+k}$ in $G_{i+k}$. Denote
this cyclic subgroup by $Z_{i+k}$.
%We will call such subgroups $Z$-\textit{subgroups}.
%Denote $G_{i,j}=\langle G_i,G_{i+1},\dots ,G_j\rangle$ for $i<j$.

(b)\qua $N=N_1\ast \dots \ast N_{k}$, where $N_l=\dots \ast
G_{l-k} \ast_{Z_l} G_l \ast_{Z_{l+k}} G_{l+k} \ast \dots $,
$(l=1,\dots ,k)$. Note that each
$N_l$ is $\psi$--invariant.

(c)\qua $N$ is the free group with the free basis $\bigcup_{i\in
\mathbb{Z}} (Y_i\setminus \{b_i\})\cup \,\{b_1,b_2,\dots ,b_k\}$.
This can be proved with the help of Tietze transformations.

For each $i\leqslant j$ denote $G_{i,j}=\langle G_i,G_{i+1},
\dots, G_j\rangle$. The group $G_{i,j}$ has two special free bases
$$\{ b_i,b_{i+1},\dots , b_{\min\{i+k-1,j\}}\}\cup \ \underset{i\leqslant l\leqslant j}{\bigcup}
(Y_l\setminus \{b_l\})$$
$$\underset{i\leqslant l\leqslant j}{\bigcup}(Y_l\setminus \{b_l\})
\cup \,\{b_j,b_{j-1},\dots ,b_{\max \{j-k+1,i\}}\}\leqno{\hbox{and}}$$ which will be
called \textit{the left} and \textit{the right basis} of $G_{i,j}$
respectively. The idea behind the left basis is the following: if
$i+k\leqslant l\leqslant j$, then we can replace each letter
$b_{l}$ by the word $b_{l-k}u_{l-k}v_{l-k}$. Thus we can
%decrease the subscript of $b_l$ and
eliminate $b_l$. The idea behind the right basis is analogous: if
$i\leqslant l\leqslant j-k$, then we can replace each letter
$b_{l}$ by the word $b_{l+k}v_l^{-1}u_l^{-1}$. In that case we can
also eliminate $b_l$.
%increase the subscript of $b_l$

Let $g$ be a nontrivial element of $N$. Consider all the subgroups
$G_{i,j}$ such that $g\in G_{i,j}$ and $j-i$ is minimal. There can
be several such subgroups (for example if $g=b_{k}b_{k+1}$, then
$g\in G_{k,k+1}$ and $g\in G_{0,1}$). Among these subgroups we
choose a subgroup with minimal $i$. Set $\alpha(g)=i$, and
$\omega(g)=j$. The number $||g||=\omega(g)-\alpha(g)+1$ will be
called the \textit{width} of the element $g$.

Denote by $g_R$ (respectively $g_L$) the cyclic reduction of $g$ written as a word in the right (in the left) basis of $G_{\alpha(g),\omega(g)}$.

%Write $g$ as a reduced word in the free group
%$G_{\alpha(g),\omega(g)}$ with respect to its left basis. Denote
%this word by $g_L$.
%%Make in this word all cyclic reductions and denote the resulting word by $g_L$.
%Analogously we can consider the right basis of
%$G_{\alpha(g),\omega(g)}$ and define the word $g_R$. For example,
%for $g=b_kb_{k+1}$ we have $g_L=g_R=u_0u_1$.

%Then $g_L$ contains a letter from $\{c_{\omega(g)},\dots
%,f_{\omega(g)} \}$. Then $g_R$ contains a letter from
%$\{c_{\alpha(g)},\dots, f_{\alpha(g)}\}$.

%Without loss of generality we assume that the words $r_i$ are
%written in the form, where $\omega(r_i)-\alpha(r_i)$ is minimal
%and $\alpha(r_i)$ is minimal. In particular any $r_i$ contains a
%letter from $Y_{\omega(r_i)}$ different from $b_{\omega(r_i)}$.
%%Not clear for $r_b>||r_i||$!!!!

\begin{defn}\label{3.1}
Let $g=z_1z_2\dots z_l$ be the normal form with respect to
the decomposition $N=N_1\ast\dots \ast N_k$, that is each $z_i$
belongs to some factor of this decomposition and $z_i,z_{i+1}$ do
not belong to the same factor. We will call any such $z_i$ \textit{a
piece of} $g$ and sometimes use $l(g)$ for $l$.

We will call $g$ a \textit{special element} of $N$ if

(1)\qua $g$ has minimal length $l$ among all its conjugates in
$N$ (this means, that $z_1$ and $z_l$ lie in different factors of
this decomposition if $l>1$),

(2)\qua if $l=1$, then $g$ has minimal width among all its
conjugates in $N$,

(3)\qua no one $z_i$ is conjugate in $H$ to a power of $b$.
\end{defn}

Note, that if $g$ is a special element and $l(g)>1$, then $g$ written
in the right basis of $G_{\alpha(g),\omega(g)}$, is cyclically
reduced. Moreover, $g$ has minimal width among all its conjugates
in $N$.

The aim of this section is to prove \fullref{3.4}. This will
be done with the help of the following lemma.

\begin{lem}\label{3.2}
Let $g$ be a special element of $N$. Then
the word $g_R$ contains a letter from
$Y_{\alpha(g)}\setminus \{b_{{\alpha(g)}}\}$ and the word $g_L$ contains a letter from
$Y_{\omega(g)}\setminus \{b_{{\omega(g)}}\}$).
\end{lem}

\begin{proof} We will prove the lemma for $g_R$.

\textbf{Case 1}\qua Suppose that $||g||\geqslant k+1$.

Then the right basis of $G_{\alpha(g),\omega(g)}$ does not contain
the letter $b_{\alpha(g)}$. Suppose that $g_R$  does not contain a
letter from $Y_{\alpha(g)}\setminus \{b_{\alpha(g)}\}$. Then
$g_R\in G_{\alpha(g)+1,\omega(g)}$, a contradiction with the
minimality of the width of $g$ among its conjugates in $N$.

\textbf{Case 2}\qua Suppose that $||g||<k+1$.

Then $G_{\alpha(g),\omega(g)}=G_{\alpha(g)}\ast \dots \ast
G_{\omega(g)}\leqslant \smash{N_{\overline{\alpha(g)}}}\ast \dots \ast
\smash{N_{\overline{\omega(g)}}}$, where $\bar{i}$ denotes the residue of
$i$ modulo $k$. Note that in this case the right basis of
$G_{\alpha(g),\omega(g)}$ coincides with $Y_{\alpha(g)}\cup \dots
\cup Y_{\omega(g)}$.
If $l(g)>1$, then as it was mentioned $g$ is cyclically reduced in this basis, and hence $g=g_R$.
 Then by condition (3), every piece of $g_R$ which lies in $\smash{G_{\alpha(g)}}$
 contains a letter from $\smash{Y_{\alpha(g)}\setminus \{b_{\alpha(g)}\}}$.
 If $l(g)=1$, then $\alpha(g)=\omega(g)$ and again by (3)
 the word $g_R$ contains a letter from $Y_{\alpha(g)}\setminus \{b_{\alpha(g)}\}$.\end{proof}

Let $B$ be a group, $A\leqslant B$, $C\vartriangleleft B $. We
will write $A\hookrightarrow B/C$ only in the case when $A\cap
C=1$, meaning the natural embedding. The following theorem is a
reformulation of the Magnus Freiheitssatz.

\begin{thm}{\rm \cite{M}}\label{3.3}\qua
 Let F be a free group with a basis $X$, and $g$ be
a cyclically reduced word in F with respect to $X$, containing a
letter $x\in X$. Then the subgroup generated by $X\setminus \{x\}$
is naturally embedded into the group $F/\langle\!\langle
g\rangle\!\rangle$.
\end{thm}

\begin{cor}\label{3.4} Let $g$ be a special element of $N$ and $j',j$ be integer numbers
such that $j'\leqslant \alpha(g)$ and $\omega(g)\leqslant j$. Then
$G_{\alpha(g)+1,j}\hookrightarrow G_{j',j}/\langle\!\langle
g\rangle\!\rangle$ and $G_{j',\omega(g)-1}\hookrightarrow
G_{j',j}/\langle\!\langle g\rangle\!\rangle$.
\end{cor}

\begin{proof}
We will prove only the first embedding. Recall that $g_R$
is the cyclic reduction of $g$ written as a word in the right basis of $G_{\alpha(g),\omega(g)}$.
The element $g_R$ remains cyclically reduced, if we rewrite it in the right basis of
$G_{j',j}$. Moreover, \fullref{3.2}
implies, that $g_R$ written in the right basis of
$G_{j',j}$ contains a letter from
$Y_{\alpha(g)}\setminus \{b_{{\alpha(g)}}\}$.
On the other hand,
any element of $G_{\alpha(g)+1,j}$ written in this basis does not
contain this letter. By \fullref{3.3},
$G_{\alpha(g)+1,j}\hookrightarrow
G_{j',j}/\langle\!\langle g\rangle\!\rangle$.
\end{proof}

\section[The structure of some quotients of G{n,m}]{The structure of some quotients of $G_{n,m}$}\label{section 4}

Let $r$ be any special element of $N$. We denote $r_i=x^{-i}rx^i$
for $i\in \mathbb{Z}$. Clearly $r_i$ is a special element.
Moreover, $\alpha(r_{i+1})=\alpha(r_i)+1$,
$\omega(r_{i+1})=\omega(r_i)+1$. In particular, all $r_i$ have the
same width. Let $j\leqslant i$. Our aim is to present the group
$G_{\alpha(r_j),\omega(r_i)}/\langle\!\langle r_j,r_{j+1},\dots
,r_i\rangle\!\rangle$ as an amalgamated product. This will be done
with the help of \fullref{4.1}.

Recall that $w_i$ denotes the word $b_i u_i v_i$ (see the notation of \fullref{section 3}).

First we introduce a technical notion: \textit{the left and the right
sets of words with respect to} $r_i$. The left set, denoted
$L(r_i)$, is $\{w_{\omega(r_i)-k},\dots ,w_{\alpha(r_i)-1}\}$. The
right set, denoted $R(r_i)$, is $\{b_{\omega(r_i)},\dots
,b_{\alpha(r_i)-1+k}\}$.
%We will denote these sets by $L(r_i)$ and $R(r_i)$ respectively.
We will assume that the subscripts of the elements of these sets
are increasing when reading from the left to the right, so these
sets are empty if $\omega(r_i)-\alpha(r_i)> k-1$. Clearly,
$L(r_i)\subset G_{-\infty,\alpha(r_i)-1}$ and $R(r_i)\subset
G_{\omega(r_i),+\infty}$.

\begin{lem}\label{4.1}
Let $r$ be a special element of $N$. Let $n,m$ and $i,j$ be
integer numbers such that $j\leqslant i$ and $m\leqslant
\alpha(r_j)$, and $\omega(r_i)\leqslant n$. Denote $s=\alpha(r_i)$
and $t=\omega(r_i)-1$. If $s>t$, we set $G_{s,t}=1$. Then the
following formula holds:
\begin{equation}\label{eq1} G_{m,n}/\langle\!\langle r_j,\dots ,r_i \rangle\!\rangle\,
\cong\, G_{m,t}/\langle\!\langle r_j,\dots
,r_{i-1}\rangle\!\rangle \mskip-25mu\underset{w_l=b_{l+k}\, (l\in L_{i,m,n})
}{{\underset{G_{s,t}}{\ast}}}\mskip-25mu G_{s,n}/\langle\!\langle
r_i\rangle\!\rangle, \end{equation} where $L_{i,m,n}=\{l\,\mid
\,w_l\in L(r_i),\, m\leqslant l\leqslant n-k\}$. Moreover, we have
\begin{equation}\label{eq2} G_{s+1,n}\hookrightarrow G_{m,n}/\langle\!\langle r_j,\dots
,r_i\rangle\!\rangle. \end{equation}
\end{lem}

\begin{proof} Note that \eqref{eq1} implies \eqref{eq2}. Indeed, by \fullref{3.4} we have the embedding $G_{s+1,n}\hookrightarrow
G_{s,n}/\langle\!\langle r_i\rangle\!\rangle$. We have the
embedding $G_{s,n}/\langle\!\langle r_i\rangle\!\rangle
\hookrightarrow G_{m,n}/\langle\!\langle r_j,\dots ,r_i
\rangle\!\rangle$ by \eqref{eq1}. Composing these two embeddings, we get the
embedding \eqref{eq2}.

Now we will prove \eqref{eq1} using induction by $i-j$. For $i-j=0$ the
formula \eqref{eq1} has the form
$$G_{m,n}/\langle\!\langle r_i\rangle\!\rangle\,
\cong\, G_{m,t}\mskip-25mu\underset{w_l=b_{l+k}\, (l\in L_{i,m,n})
}{{\underset{G_{s,t}}{\ast}}}\mskip-25mu G_{s,n}/\langle\!\langle
r_i\rangle\!\rangle.$$ Let $M$ be the subgroup of $G_{m,n}$
generated by $G_{s,t}$ and the set $\{ b_{l+k}\,|\, l\in
L_{i,m,n}\}$. Clearly, $M$ is a subgroup of $G_{m,t}$ and
$G_{s,n}$. It is sufficient to prove that $M$ embeds into
$G_{s,n}/\langle\!\langle r_i\rangle\!\rangle$. Consider two
cases.

\textbf{Case 1}\qua Suppose that $k\geqslant ||r_i||$.

By definition, the group $M$ lies in the subgroup generated by the
set $Y_{\alpha(r_i)}\cup\dots \cup Y_{\omega(r_i)-1}\cup
\,\{b_{\alpha(r_i)},\dots, b_{\min \{\alpha(r_i)-1+k,n\}}\}$. In
that case this set is a part of the left basis of
$G_{s,n}$. Consider the cyclically reduced word in this basis,
corresponding to $r_i$. \fullref{3.2} implies that it contains a letter
from $Y_{\omega(r_i)}\setminus \{b_{\omega(r_i)}\}$. Hence, by
\mbox{\fullref{3.3}}, $M$ embeds into $G_{s,n}/\langle\!\langle
r_i\rangle\!\rangle$.

\textbf{Case 2}\qua Suppose that $k<||r_i||$.

In that case $M=G_{s,t}$ and the desired embedding follows from
\fullref{3.4}.

Thus, the base of induction holds. Suppose that the formula \eqref{eq1}
holds for $j,i$ and prove it for $j,i+1$. Thus, we need to prove
that
\begin{equation}\label{eq3} G_{m,n}/\langle\!\langle r_j,\dots ,r_{i+1} \rangle\!\rangle\,
\cong\, G_{m,t+1}/\langle\!\langle r_j,\dots ,r_i\rangle\!\rangle
\mskip-25mu\underset{w_l=b_{l+k}\, (l\in L_{i+1,m,n})
}{{\underset{G_{s+1,t+1}}{\ast}}}\mskip-25mu G_{s+1,n}/\langle\!\langle
r_{i+1}\rangle\!\rangle.\end{equation}
Let $M$ be the subgroup of $G_{m,n}$ generated by $G_{s+1,t+1}$
and the set $\{ w_l\,|\, l\in L_{i+1,m,n}\}$. Equivalently, $M$ is generated by $G_{s+1,t+1}$ and the set
$\{ b_{l+k}\,|\, l\in L_{i+1,m,n}\}$. It is sufficient to prove
that $M$ embeds naturally into the factors of \eqref{eq3}, that is into
$G_{m,t+1}/\langle\!\langle r_j,\dots ,r_i\rangle\!\rangle$ and
$G_{s+1,n}/\langle\!\langle r_{i+1}\rangle\!\rangle$.

%By Claim 1 we have
%$$M=G_{s+1,t+1}\ast F,$$ where $F$ is the free group with the free
%basis $\{ u_l\,|\, l\in L_{i+1,n}\}$.
The group $M$ can be considered as a subgroup of $G_{s+1,n}$, and
$G_{s+1,n}$ naturally embeds into $G_{m,n}/\langle\!\langle
r_j,\dots ,r_i\rangle\!\rangle$ by \eqref{eq2}. Hence $M$ naturally embeds
into $G_{m,n}/\langle\!\langle r_j,\dots ,r_i\rangle\!\rangle$.
%Also $M$ can be considered as a subgroup of $G_{m,t+1}$.
Thus $M$
naturally embeds into $G_{m,t+1}/\langle\!\langle r_j,\dots
,r_i\rangle\!\rangle$, since $M\leqslant G_{m,t+1}\leqslant G_{m,n}$.

The embedding of $M$ into $G_{s+1,n}/\langle\!\langle
r_{i+1}\rangle\!\rangle$ can be proved by the same argument as in
the case of the base of induction. \end{proof}

The following lemma can be proved similarly.

\begin{lem}\label{4.2} Let $r$ be a special element of $N$.
Let $n,m$ and $i,j$ be integer numbers such that $j\leqslant i$
and $m\leqslant \alpha(r_j)$, and $\omega(r_i)\leqslant n$. Denote
$s=\alpha(r_j)+1$ and $t=\omega(r_i)$. If $s>t$, we set
$G_{s,t}=1$. Then the following formula holds:
$$G_{m,n}/\langle\!\langle r_j,\dots ,r_i \rangle\!\rangle\,
\cong\, G_{m,t}/\langle\!\langle r_j\rangle\!\rangle\mskip-25mu
\underset{w_l=b_{l+k}\, (l\in L_{j,m,n})
}{{\underset{G_{s,t}}{\ast}}}\mskip-25mu G_{s,n}/\langle\!\langle
r_{j+1},\dots ,r_i\rangle\!\rangle , $$ where
$L_{i,m,n}=\{l\,\mid \,w_l\in L(r_i),\, m\leqslant l\leqslant
n-k\}$.
\end{lem}

\section[Proof of Proposition 2.1]{Proof of \fullref{2.1}}\label{section 5}

Let $r$ and $s$ be two elements of $H$ with the same normal
closure and $r_x=0$. Recall that $N$ denotes the kernel of the
homomorphism $H\rightarrow \mathbb{Z}$, sending $x$ to 1 and each
other generator of $H$ to 0.  Denote $r_i=x^{-i}rx^i$,
$s_i=x^{-i}sx^i$, $i\in \mathbb{Z}$. Then $r_i,s_i\in N$.
Moreover, the sets $\mathcal{R}=\{\dots , r_{-1},r_0,r_1,\dots \}$
and $\mathcal{S}=\{ \dots ,s_{-1},s_0,s_1,\dots \}$ have the same
normal closure in $N$. We will prove that some $r_i$ is conjugate
to $s_0^{\pm 1}$ in $N$. This will imply, that  $r$ is conjugate
to $s^{\pm 1}$ in $H$.

We may assume, that $r$ and $s$ are special elements. Indeed, let
$r=z_1z_2\dots z_l$ and $s=c_1c_2\dots c_{l'}$ be normal forms
with respect to the decomposition $N_1\ast \dots \ast N_k$.
Conjugating, we may assume that the condition (1) of \fullref{3.1} is satisfied. Applying a power of an automorphism $\psi$
from \fullref{2.2}, we may assume that the condition (3) is
satisfied. Finally, if $l=1$ or $l'=1$, we may conjugate $r$ or
$s$ to ensure the condition (2).

It follows that $r_i$ and $s_i$ are special elements and
$\alpha(r_{i+1})=\alpha(r_i)+1$, $\omega(r_{i+1})=\omega(r_i)+1$.
In particular, all $r_i$ have the same width. The same is valid
for $s_i$.

Since $s_0$ can be deduced from $\mathcal{R}$ in $N$, there exist
integer numbers $j,i$ such that $j\leqslant i$ and $s_0$ is
trivial in $G_{\alpha,\omega}/\langle\!\langle r_j,r_{j+1},\dots
,r_i \rangle\!\rangle$, where $\alpha=\alpha(r_j)$,
$\omega=\omega(r_i)$. We assume that $i-j$ is minimal possible. By
\fullref{4.1} we have $$G_{\alpha,\omega}/\langle\!\langle
r_j,\dots ,r_i \rangle\!\rangle\, \cong\,
G_{\alpha,\omega-1}/\langle\!\langle r_j,\dots
,r_{i-1}\rangle\!\rangle \underset{A}{\ast}
G_{\alpha(r_i),\omega}/\langle\!\langle r_i\rangle\!\rangle,$$ for
some subgroup $A$.

It follows that $s_0\notin G_{\alpha,\omega-1}$, otherwise $s$
were trivial in $G_{\alpha,\omega-1}/\langle\!\langle r_j,\dots
,r_{i-1}\rangle\!\rangle$, that contradicts to the minimality of
$i-j$. Hence, $s_0$ written as a word in the left basis of
$G_{\alpha,\omega}$ must contain a letter from $Y_{\omega}$.

Now we will prove that $\alpha(s_0)=\alpha$ and
$\omega(s_0)=\omega$. If $s_0$ contains a letter from
$Y_{\omega}\setminus \{b_{\omega}\}$, then clearly,
$\omega(s_0)\geqslant \omega$.

Suppose that $s_0$ contains the letter $b_{\omega}$, but does not
contain any letter from $Y_{\omega}\setminus \{b_{\omega}\}$.
%We can replace $b_{\omega}$ by $w_{\omega}$.
Then $b_{\omega}$ belongs to the left basis of
$G_{\alpha,\omega}$, what can happens only if
$\omega-\alpha+1\leqslant k$. But in this case $s_0$ contains a
piece, which is a power of $b_{\omega}$ -- a contradiction.

%If $\omega-\alpha+1> k$, then $b_{\omega}=w_{\omega-k}\in
%G_{\alpha,\omega-1}$ and hence $s_0\in G_{\alpha,\omega-1}$ -- a
%contradiction.

%If $\omega-\alpha+1\leqslant k$, then $s_0$ contains a piece,
%which is a power of $b_{\omega}$ -- a contradiction.

Thus, we have proved that $\omega(s_0)\geqslant\omega $.
Analogously, $\alpha(s_0)\leqslant \alpha$. Hence,
$\alpha(s_0)=\alpha$ and $\omega(s_0)=\omega $. In particular,
$||s_0||=\omega-\alpha+1\geqslant ||r_j||$. By symmetry,
$||r_j||\geqslant ||s_0||$. Hence $||r_j||=||s_0||$ and
$\alpha=\alpha(r_j)$, $\omega=\omega(r_j)$. It follows that $s_0$
can be deduced from $r_j$ in $G_{\alpha,\omega}$ and the subscript
$j$ is determined from the equation $\alpha(s_0)=\alpha(r_j)$.
Similarly, $r_j$ can be deduced in $G_{\alpha,\omega}$ from $s_0$.
By \fullref{free}, $s_0$ is conjugate to $r_j^{\pm 1}$ in
$G_{\alpha,\omega}$. Hence $s$ is conjugate to $r^{\pm 1}$ in $H$.
\qed

\bibliographystyle{gtart}
\bibliography{link}

\end{document}